\documentclass[12pt,a4paper]{article}

\usepackage[margin=1in]{geometry}  
\usepackage{graphicx}              
\usepackage{amsmath}               
\usepackage{amsfonts}              
\usepackage{amsthm}                
\usepackage{amssymb}

\usepackage{mathrsfs}
\usepackage{dsfont}

\usepackage{hyperref}

\usepackage{authblk}

\usepackage{verbatim}

\usepackage[style=authoryear]{biblatex}

\DeclareMathOperator{\Tr}{Tr}
\DeclareMathOperator{\Rk}{Rk}
\newcommand{\bd}[1]{\mathbf{#1}}   

\renewcommand{\baselinestretch}{1.2}

\begin{document}

\title{A note on random samples of Lie algebras}

\author{Georg Beyerle\thanks{Address correspondence to: mail@gbeyerle.de}}

\affil{Caputh, Germany}

\date{7 July 2014}

\maketitle

\begin{abstract}
Recently, Paiva and Teixeira (arXiv:1108.4396)
showed that the structure constants of a Lie algebra
are the solution of a system of linear equations
provided a certain subset of the structure constants
are given a-priori.
Here it is noted
that Lie algebras generated in this way
are solvable and
their derived subalgebras are Abelian
if the system of linear equations
considered by Paiva and Teixeira is not degenerate.
An efficient numerical algorithm
for the calculation of their structure constants
is described.
\end{abstract}

\section{Introduction}
\label{sect:intro}

Recently,
\textcite{paiva:2011}
showed that the structure constants
of a finite-dimensio\-nal Lie algebra
are the solution of a system of linear equations,
provided a certain subset of structure constants
are given a-priori.
We recall that a Lie algebra~$\mathfrak{L}$
is a vector space over a field~$\mathbb{F}$
equipped with a bilinear product,
the Lie bracket $[\;x\;,\;y\;]$ with $x,y \in \mathfrak{L}$,
\parencite[see, e.g.,][]{hall:2003,%
humphreys:1972,knapp:2005,jacobson:1979}.
The following discussion is restricted to finite-dimensional
Lie algebras over the field of
the real ($\mathbb{F} = \mathbb{R}$) or
the complex numbers ($\mathbb{F} = \mathbb{C}$).
The Lie bracket satisfies
$[\;x\;,\;x\;] = 0$ for all $x \in \mathfrak{L}$
and the Jacobi identity
\begin{eqnarray}
\label{eq:jacobiident}
 [ \; x \; , \; [ \; y \; , \; z \; ] \; ] +
 [ \; z \; , \; [ \; x \; , \; y \; ] \; ] +
 [ \; y \; , \;  [\; z \; , \; x \; ] \; ]
 &=& 0
\end{eqnarray}
for all $x,y,z \in \mathfrak{L}$.
The $N$-dimensional Lie algebra~$\mathfrak{L}$
is completely characterized by
the coordinates (known as structure constants)~$\bd{f}\{i,j,k\}$
of the Lie bracket product
with respect to the basis~$g_i$
\parencite[see, e.g.,][]{hall:2003,humphreys:1972,jacobson:1979},
\begin{eqnarray}
\label{eq:liestrconst}
  [ \; g_i \; , \; g_j \; ] = \sum_{k=1}^N \bd{f}\{i,j,k\} \; g_k
  \qquad i,j=1,\ldots,N
  \quad.
\end{eqnarray}
Here and in the following, the notation
$\bd{X}\{i_1,\ldots,i_D\}$
is used to identify a specific element
of the $D$-dimensional matrix~$\bd{X}$
indexed by $D$ positive integers $i_1,\ldots,i_D$.
Vectors are represented as $2$-dimensional matrices
with one singleton dimension.
Specifically,
$i_1=1$ and $i_2=1,\ldots,N$
describes a row vector
$\vec{\bd{X}} \equiv
\left(
  \bd{X}\{1,1\},\ldots,\bd{X}\{1,N\}
\right)$
and
$i_1=1,\ldots,N$ and $i_2=1$
its transpose,
a column vector
$\vec{\bd{X}}^{\;T} \equiv
\left(
  \bd{X}\{1,1\},\ldots,\bd{X}\{N,1\}
\right)$.

\textcite{paiva:2011}
observed that the Jacobi identity
(\ref{eq:jacobiident}),
expressed in terms
of the structure constants~$\bd{f}\{i,j,k\}$, is
\begin{eqnarray}
\label{eq:jacobistructconst}
  J^A_{i,j,k,m} + J^B_{i,j,k,m}  + J^C_{i,j,k,m}  &=& 0
\end{eqnarray}
with
\begin{eqnarray}
\label{eq:definejsums}
J^A_{i,j,k,m}
  \equiv
  \sum_{l=1}^N \bd{f}\{i,j,l\} \, \bd{f}\{k,l,m\}
  \\
J^B_{i,j,k,m}
  \equiv
  \sum_{l=1}^N \bd{f}\{k,i,l\} \, \bd{f}\{j,l,m\}
  \nonumber \\
J^C_{i,j,k,m}
  \equiv
  \sum_{l=1}^N \bd{f}\{j,k,l\} \, \bd{f}\{i,l,m\}
  \nonumber
\end{eqnarray}
and
$i,j,k,m=1,\ldots,N$,
represents a system of linear equations
in unknowns
\begin{eqnarray}
\label{eq:unknownvars}
\bd{f}\{2\le i\le N,i+1\le j\le N,1\le k\le N\}
\end{eqnarray}
if the index~$i$ in (\ref{eq:jacobistructconst})
is fixed at $i=1$.
Since
$\bd{f}\{i,j,k\} = -\bd{f}\{j,i,k\}$,
the three summands in
(\ref{eq:jacobistructconst})
are related,
\begin{eqnarray}
\label{eq:jantisym}
J^A_{i,j,k,m}
  &=& -J^A_{j,i,k,m} = -J^C_{k,j,i,m} = -J^B_{i,k,j,m}
  \\
J^B_{i,j,k,m}
  &=& -J^C_{j,i,k,m} = -J^B_{k,j,i,m} = -J^A_{i,k,j,m}
  \nonumber \\
J^C_{i,j,k,m}
  &=& -J^B_{j,i,k,m} = -J^A_{k,j,i,m} = -J^C_{i,k,j,m}
  \nonumber
\end{eqnarray}
and the range of $i$, $j$ and $k$ in (\ref{eq:jacobistructconst})
can be restricted to $1 \le i < j < k \le N$
without omitting linear independent equations.
\textcite{paiva:2011}
considered the system of equations
\begin{eqnarray}
\label{eq:syslineqn}
  J^A_{i=1,j,k,m} + J^B_{i=1,j,k,m}  + J^C_{i=1,j,k,m}
  &=& \sum_{l=1}^N \;\bd{f}\{i=1,j,l\} \, \bd{f}\{k,l,m\}
  \\
  && - \; \bd{f}\{i=1,k,l\} \, \bd{f}\{j,l,m\}
  + \bd{f}\{i=1,l,m\} \, \bd{f}\{j,k,l\}
  \nonumber \\
  &=& 0
 \nonumber
\end{eqnarray}
with
\begin{eqnarray}
\label{eq:syslineqn2}
  && j=2,\ldots,N
  \qquad
  k=j+1,\ldots,N
  \qquad
  m=1,\ldots,N
\end{eqnarray}
and
$\bd{f}\{i=1,j,k\}$
taken as known parameters.

The number of individual linear equations
in (\ref{eq:syslineqn}),
\begin{eqnarray}
\label{eq:nofeqn}
N \, \sum_{j=2}^N \left( \sum_{k=j+1}^N 1 \right)
  = N \, \sum_{j=2}^N \left(N - j\right)
  = N\,(N-1)\,(N-2)\,/\,2
\end{eqnarray}
matches the number of unknowns (\ref{eq:unknownvars}).
Thus, a unique solution exists
provided the a-priori parameters~$\bd{f}\{i=1,j,k\}$
are chosen such,
that the system of linear equations is non-degenerate
\parencite[see, e.g,][]{strang:2009}.
We assume that this condition is fulfilled.

In the following section
an efficient algorithm for calculating
$\bd{f}\{i,j,k\}$ is described
bypassing the task of solving
the system of linear equations
in (\ref{eq:syslineqn}).
The appendix includes
a computer implementation of the algorithm.
Second, it is shown
that Lie algebras randomly generated in this way
belong to the class of solvable Lie algebras
and their derived subalgebras
$[\;\mathfrak{L}\;,\;\mathfrak{L}\;]$
are Abelian.
The algorithm described by
\textcite{paiva:2011} may therefore be used
to generate random solvable Lie algebras
\parencite[as regards random samples of
nilpotent Lie algebras, see also][]{luks:1977}.

\section{Generating a random Lie algebra}
\label{sect:algorithm}

As will be evident in the following,
it is convenient to rewrite
the problem of solving~(\ref{eq:syslineqn}),
in terms of the $N \times N$ matrices $\bd{A}_k$
with matrix elements
\begin{eqnarray}
\label{eq:adjrepvsstrconst}
\bd{A}_k\{i,j\}
  &=& \bd{f}\{k,j,i\} \quad .
\end{eqnarray}
$\bd{A}_k$ ($k=1,\ldots,N$)
are the adjoint representation of~$\mathfrak{L}$,
\parencite[see, e.g.,][]{hall:2003,humphreys:1972}
and
(\ref{eq:liestrconst})
translates into the matrix equation
\begin{eqnarray}
\label{eq:adjrepdefeqn}
[\; \bd{A}_i \;,\; \bd{A}_j \;] =
\bd{A}_i \cdot \bd{A}_j - \bd{A}_j \cdot \bd{A}_i
&=& \sum_{k=1}^N \bd{f}\{i,j,k\} \; \bd{A}_k
\\
&=& \sum_{k=1}^N \bd{A}_i\{k,j\} \; \bd{A}_k
\nonumber
\quad.
\end{eqnarray}
In terms of the adjoint representation
the a-priori parameters
are precisely the
elements of the matrix, $\bd{A}_{k=1}$.
According to
\textcite{paiva:2011}
all matrices of the adjoint representation
can be calculated from $\bd{A}_{k=1}$
provided the constraints
discussed above are satisfied.

We introduce an a-priori matrix~$\bd{P}$
containing $N\,(N-1)$ parameters
\begin{eqnarray}
\label{eq:apriorip}
  c \; \bd{A}_1
  \equiv
  \bd{P}
  \equiv
  \left(
    \begin{array}{ccccc}
    0      & & \bd{P}\{1,2\} & \ldots & \bd{P}\{1,N\}\\
    \vdots & & \vdots        &        & \vdots \\
    0      & & \bd{P}\{N,2\} & \ldots & \bd{P}\{N,N\}\\
    \end{array}
  \right)
\end{eqnarray}
with some scalar factor~$c \in \mathbb{F}$.
Since $\bd{f}\{i,i,k\} = 0$,
the first column
of $\bd{A}_{k=1}$ is zero and,
therefore,
the rank of $\bd{P}$, $\Rk(\bd{P})$, is at most $N-1$;
in the following it is assumed
that the parameters $\bd{P}\{i,j\}$ are chosen such
that
\begin{eqnarray}
\label{eq:rankcond}
\Rk(\bd{P}) = N-1
\quad.
\end{eqnarray}
The null space of $\bd{P}$
is therefore one-dimensional
and $\vec{\bd{n}}$,
the normalized null vector of $\bd{P}$,
is unique (up to its sign) with
\begin{eqnarray}
\label{eq:nulldotp}
\vec{\bd{n}} \cdot \bd{P} = (0,\ldots,0)
\quad.
\end{eqnarray}
Following \textcite{paiva:2011}
we expect the elements of the matrices $\bd{A}_{k\ne1}$
to solve a system of linear equations
and introduce the {\sl Ansatz}
\begin{eqnarray}
\label{eq:adjmatfromt}
\bd{A}_k
  &=& \bd{P} \cdot \bd{T}_k \qquad k=1,\ldots\,N
\end{eqnarray}
with $\bd{T}_k$ given by
\begin{eqnarray}
\label{eq:tmatansatz}
\bd{T}_k
  = \bd{n}\{k\} \; \mathds{1} -
    \vec{\bd{e}}_k^{\;T} \cdot \vec{\bd{n}}
\quad.
\end{eqnarray}
Here,
$\mathds{1}$ is the $N\times N$ unit matrix
and $\vec{\bd{e}}_k$ is the row vector with $N$ elements,
\begin{eqnarray}
\label{eq:defvectek}
\vec{\bd{e}}_k
  &=& (0,\ldots,0,\underset{k}{1},0,\ldots,0)
\end{eqnarray}
I.e.\ $\bd{e}_k\{k\}=1$ and $\bd{e}_k\{i\}=0$ for all $i\ne k$.
Note that
$\vec{\bd{n}}\cdot\vec{\bd{e}}_k^{\;T}$
is a scalar, whereas
$\vec{\bd{e}}_k^{\;T} \cdot \vec{\bd{n}}$
evaluates to an $N\times N$ matrix.

It is convenient to choose
\begin{eqnarray}
\label{eq:scalfactorc}
c \equiv \frac{1}{\bd{n}\{1\}}
\end{eqnarray}
in (\ref{eq:apriorip})
assuming $\bd{n}\{1\} \ne 0$.
The choice (\ref{eq:scalfactorc}) allows us to write
both, $\bd{A}_{k=1}$ and $\bd{A}_{k\ne1}$
in the form of (\ref{eq:adjmatfromt}),
since
\begin{eqnarray}
\label{eq:lcomm1}
  \bd{A}_{k=1}
  &=&
  \bd{P} \cdot \bd{T}_{k=1} \\
  &=&
  \bd{P} \cdot
  \left(
    \begin{array}{ccccc}
    0      & -\bd{n}\{2\} &    \ldots   &    \ldots      & -\bd{n}\{N\} \\
    0      &  \bd{n}\{1\} &       0     &    \ldots      &       0      \\
    \vdots &      0       &    \ddots   &                &    \vdots    \\
    \vdots &   \vdots     &             &    \ddots      &       0      \\
    0      &      0       &    \ldots   &      0         & \bd{n}\{1\}  \\
    \end{array}
  \right)
  \nonumber \\
  &=&
  \bd{n}\{1\} \; \bd{P}
  \nonumber
\end{eqnarray}
and taking into account
the vanishing first column of~$\bd{P}$.
In order to prove that {\sl Ansatz}
(\ref{eq:tmatansatz}) indeed constitutes
an adjoint representation of~$\mathfrak{L}$,
the matrices $\bd{A}_k = \bd{P}\cdot\bd{T}_k$
must be shown to satisfy (\ref{eq:adjrepdefeqn}).
We proceed by first calculating
the left-hand side of (\ref{eq:adjrepdefeqn}),
the Lie bracket of $\bd{A}_i$ and $\bd{A}_j$,
\begin{eqnarray}
\label{eq:lcomm1}
[\;\bd{A}_i\;,\;\bd{A}_j\;]
  &=& \bd{A}_i \cdot \bd{A}_j -
      \bd{A}_j \cdot \bd{A}_i \\
  &=&  \quad \bd{P} \cdot
         \left( \bd{n}\{i\} \; \mathds{1} -
                \vec{\bd{e}}_i^{\;T} \cdot \vec{\bd{n}}
         \right) \,
       \cdot \bd{P} \cdot
         \left( \bd{n}\{j\} \; \mathds{1} -
                \vec{\bd{e}}_j^{\;T} \cdot \vec{\bd{n}}
         \right) \nonumber \\
  && - \; \bd{P} \cdot
       \left( \bd{n}\{j\} \; \mathds{1}
         - \vec{\bd{e}}_j^{\;T} \cdot \vec{\bd{n}}
       \right)
       \cdot \bd{P} \cdot
       \left(\bd{n}\{i\} \; \mathds{1}
         - \vec{\bd{e}}_i^{\;T} \cdot \vec{\bd{n}}
       \right) \nonumber \\
  &=&  \quad \bd{n}\{i\} \,
       \bd{P}^2 \cdot
       \left( \bd{n}\{j\} \; \mathds{1}
         - \vec{\bd{e}}_j^{\;T} \cdot \vec{\bd{n}}
       \right) \nonumber \\
  && - \; \bd{n}\{j\} \,
       \bd{P}^2 \cdot
       \left(\bd{n}\{i\} \; \mathds{1}
         - \vec{\bd{e}}_i^{\;T} \cdot \vec{\bd{n}}
       \right) \nonumber
\end{eqnarray}
on account of (\ref{eq:nulldotp}) and, finally,
\begin{eqnarray}
\label{eq:lcomm2}
[\;\bd{A}_i\;,\;\bd{A}_j\;]
  &=&  \quad \bd{n}\{i\} \, \bd{n}\{j\} \, \bd{P}^2
     - \bd{n}\{i\} \, \bd{P}^2 \cdot \vec{\bd{e}}_j^{\;T}\cdot\vec{\bd{n}}
  \\
  && - \; \bd{n}\{j\} \, \bd{n}\{i\} \, \bd{P}^2
     + \bd{n}\{j\} \, \bd{P}^2 \cdot \vec{\bd{e}}_i^{\;T}\cdot\vec{\bd{n}}
  \nonumber \\
  &=& \bd{P}^2 \cdot
      \left( \bd{n}\{j\} \; \vec{\bd{e}}_i^{\;T} -
             \bd{n}\{i\} \; \vec{\bd{e}}_j^{\;T} \right)
      \cdot \vec{\bd{n}} \nonumber
\quad.
\end{eqnarray}
Second,
the right-hand side of (\ref{eq:adjrepdefeqn})
is found to be
\begin{eqnarray}
\label{eq:ansatzproof}
\sum_{k=1}^N \bd{A}_i\{k,j\} \, \bd{A}_k
&=& \sum_{k=1}^N
  \left(
    \bd{P}\,(\bd{n}\{i\} \; \mathds{1} -
    \vec{\bd{e}}_i^{\;T}\cdot\vec{\bd{n}})
  \right)\{k,j\} \;\;
  \bd{P} \cdot
  \left(
    \bd{n}\{k\} \; \mathds{1} -
                    \vec{\bd{e}}_k^{\;T}\cdot\vec{\bd{n}}
  \right)
  \\
&=& \sum_{k=1}^N
    \left(
      \bd{P}\,(\bd{n}\{i\} \; \mathds{1} -
      \vec{\bd{e}}_i^{\;T}
      \cdot\vec{\bd{n}})\cdot
      \vec{\bd{e}}_j^{\;T}
    \right)\{k\} \;\;
    \bd{P} \cdot
    \left(
      \bd{n}\{k\} \; \mathds{1} -
                 \vec{\bd{e}}_k^{\;T}\cdot\vec{\bd{n}}
      \right)
    \nonumber \\
&=& \quad \bd{P} \cdot \sum_{k=1}^N
    \left(
      \bd{P}\,(\bd{n}\{i\} \; \mathds{1} -
      \vec{\bd{e}}_i^{\;T}
      \cdot\vec{\bd{n}})\cdot\vec{\bd{e}}_j^{\;T}
    \right)\{k\} \;\;
    \bd{n}\{k\}
    \nonumber \\
&& - \; \bd{P} \cdot \sum_{k=1}^N
    \left(
      \bd{P}\,(\bd{n}\{i\} \; \mathds{1}-
      \vec{\bd{e}}_i^{\;T}
      \cdot\vec{\bd{n}})\cdot\vec{\bd{e}}_j^{\;T}
    \right)\{k\} \;\;
    \vec{\bd{e}}_k^{\;T} \cdot \vec{\bd{n}}
    \nonumber \\
&=& \quad \bd{P} \cdot \vec{\bd{n}} \cdot
    \bd{P}\,(\bd{n}\{i\} \; \mathds{1} - \vec{\bd{e}}_i^{\;T}
    \cdot \vec{\bd{n}}) \cdot \vec{\bd{e}}_j^{\;T}
    \nonumber \\
&& - \; \bd{P} \cdot
    \left(\sum_{k=1}^N
      \left( \bd{P} \cdot (\bd{n}\{i\} \; \mathds{1} -
        \vec{\bd{e}}_i^{\;T}
        \cdot\vec{\bd{n}}) \cdot \vec{\bd{e}}_j^{\;T}
      \right)\{k\} \;\;
      \vec{\bd{e}}_k^{\;T}
    \right) \cdot \vec{\bd{n}}
    \nonumber
\quad.
\end{eqnarray}
The first term vanishes
owing to (\ref{eq:nulldotp})
and the second simplifies with (\ref{eq:defvectek}) to
\begin{eqnarray}
\sum_{k=1}^N \bd{A}_i\{k,j\} \, \bd{A}_k
&=& -\bd{P} \cdot
    \left(
      \bd{P} \cdot (\bd{n}\{i\} \; \mathds{1} - \vec{\bd{e}}_i^{\;T}
      \cdot \vec{\bd{n}}) \cdot \vec{\bd{e}}_j^{\;T}
    \right) \cdot \vec{\bd{n}}
    \\
&=& \bd{P}^2 \cdot
    \left(
      \bd{n}\{j\}\,\vec{\bd{e}}_i^{\;T} -
      \bd{n}\{i\}\,\vec{\bd{e}}_j^{\;T}
    \right)
    \cdot \vec{\bd{n}}
    \nonumber
\end{eqnarray}
which equals the right-hand side of (\ref{eq:lcomm2}).
We conclude that the matrices
\begin{eqnarray}
\label{eq:adjmatfromt2}
\bd{A}_k
  &=& \bd{P} \cdot
      \left(
        \bd{n}\{k\} \; \mathds{1} -
        \vec{\bd{e}}_k^{\;T} \cdot \vec{\bd{n}}
      \right) \qquad k=1,\ldots\,N
\end{eqnarray}
form an adjoint representation of~$\mathfrak{L}$.
Furthermore, the proof of {\sl Ansatz} (\ref{eq:adjmatfromt})
implies that the adjoint representation matrices of~$\mathfrak{L}$
(and with the help of (\ref{eq:adjrepvsstrconst})
also all its structure constants)
can be uniquely determined,
once the matrix~$\bd{P}$
with its $N\,(N-1)$ a-priori parameters is given,
provided
$\Rk(\bd{P}) = N-1$
and $\bd{n}\{1\} \ne 0$.

\section{Derived series and lower central series}
\label{sect:solvability}

Knowing the matrices of the adjoint representation
(\ref{eq:adjmatfromt2})
it is straightforward to calculate
the derived series of $\mathfrak{L}$
\parencite{humphreys:1972}.
With the abbreviation
\begin{eqnarray}
\label{eq:abbr}
\vec{\bd{m}}_{j,i}^T
  \equiv \bd{n}\{j\} \cdot \vec{\bd{e}}_i^{\;T}
       - \bd{n}\{i\} \cdot \vec{\bd{e}}_j^{\;T}
\end{eqnarray}
from (\ref{eq:lcomm2}) it follows
\begin{eqnarray}
\label{eq:derivser}
[ \; [ \; \bd{A}_i\;,\;\bd{A}_j \; ] \; ,
       \; [ \; \bd{A}_k \; , \; \bd{A}_l \; ] \; ]
  &=& [ \; \bd{P}^2 \cdot
      \vec{\bd{m}}_{j,i}^T \cdot \vec{\bd{n}} \; , \;
         \bd{P}^2 \cdot
      \vec{\bd{m}}_{l,k}^T \cdot \vec{\bd{n}} \; ] \\
  &=& \quad\bd{P}^2 \cdot
      \vec{\bd{m}}_{j,i}^T \cdot \vec{\bd{n}} \cdot
         \bd{P}^2 \cdot
      \vec{\bd{m}}_{l,k}^T \cdot \vec{\bd{n}} \nonumber \\
  &&  - \;\bd{P}^2 \cdot
      \vec{\bd{m}}_{l,k}^T \cdot \vec{\bd{n}} \cdot
         \bd{P}^2 \cdot
      \vec{\bd{m}}_{j,i}^T \cdot \vec{\bd{n}} \nonumber \\
  &=& 0 \nonumber
\end{eqnarray}
for all $i,j,k,l=1,\ldots,N$
owing to (\ref{eq:nulldotp}).
Thus, the Lie algebra~$\mathfrak{L}$ is solvable
and the derived subalgebra~$[\;\mathfrak{L}\;,\;\mathfrak{L}\;]$
is Abelian.

Alternatively, solvability of $\mathfrak{L}$
may be proved using  Cartan's criterion,
i.e.\ the fact that
the Killing form~$K(x,y)$ with
$x \in \mathfrak{L}$ and
$y \in [\;\mathfrak{L}\;,\;\mathfrak{L}\;]$
is identically zero
if and only if $\mathfrak{L}$ is solvable
\parencite{knapp:2005}.
We find
\begin{eqnarray}
\label{eq:killfrm}
\Tr\left(
     \bd{A_i} \cdot [ \; \bd{A}_j \; , \; \bd{A}_k \; ]
   \right)
  &=&
  \Tr\left(
       [ \; \bd{A}_j \; , \; \bd{A}_k \; ] \cdot \bd{A_i}
     \right)
  \\
  &=&
  \Tr \left(
        \bd{P}^2 \cdot
        \vec{\bd{m}}_{k,j}^T \cdot \vec{\bd{n}}
        \cdot
        \bd{P} \cdot
        \left(
          \bd{n}\{i\} \; \mathds{1} -
          \vec{\bd{e}}_i^{\;T} \cdot \vec{\bd{n}}
        \right)
      \right)
  \nonumber \\
  &=& \Tr(0) = 0
  \nonumber
\end{eqnarray}
for all $i,j,k=1,\ldots,N$
owing to (\ref{eq:nulldotp}).

The lower central series
\parencite{humphreys:1972},
on the other hand, is
\begin{eqnarray}
\label{eq:lowcentral2}
[\;\bd{A}_{i_L}\;,\;\ldots\;[\;\bd{A}_{i_1}\;,
  \;[\;\bd{A}_j\;,\;\bd{A}_k\;]\;]\;\ldots\;]
  &=& \bd{n}\{i_L\} \ldots \bd{n}\{i_1\} \; \bd{P}^{L+2} \cdot
      \vec{\bd{m}}_{k,j}^T \cdot \vec{\bd{n}}
\end{eqnarray}
which follows from
\begin{eqnarray}
\label{eq:lowcentral1}
[\;\bd{A}_i\;,\;[\;\bd{A}_j\;,\;\bd{A}_k\;]\;]
  &=& [ \; \bd{P} \cdot
    ( \bd{n}\{i\} \; \mathds{1} -
      \vec{\bd{e}}^{\;T}_i\cdot\vec{\bd{n}} ) \; , \;
         \bd{P}^2 \cdot
      \vec{\bd{m}}_{k,j}^T \cdot \vec{\bd{n}}\;]  \\
  &=& \bd{n}\{i\} \; \bd{P}^3 \cdot
      \vec{\bd{m}}_{k,j}^T \cdot \vec{\bd{n}}
      \nonumber
\quad.
\end{eqnarray}
If the parameter matrix~$\bd{P}$ is not nilpotent,
the right-hand side of (\ref{eq:lowcentral2})
does not vanish.
Accordingly, the lower central series of $\mathfrak{L}$
will not terminate
and $\mathfrak{L}$ not be nilpotent.

If, however, $\bd{P}$
is chosen to be strictly upper-triangular
(upper-triangular with zeros on the main diagonal),
then $\bd{P}^N = 0$.
In this case the series (\ref{eq:lowcentral2})
terminates and
the corresponding Lie algebra is nilpotent, too.

\section{Concluding remarks}
\label{sect:discussion}

We have shown that
the adjoint representation of a
random, solvable Lie algebra~$\mathfrak{L}$
can be obtained from a (real or complex)
parameter matrix $\bd{P}$
\begin{eqnarray}
\label{eq:adjmatfromt3}
\bd{A}_k
  &=& \bd{P} \cdot
      \left(
        \bd{n}\{k\} \; \mathds{1} -
        \vec{\bd{e}}_k^{\;T} \cdot \vec{\bd{n}}
      \right) \qquad k=1,\ldots\,N
\end{eqnarray}
where $\vec{\bd{n}}$ denotes the null vector of $\bd{P}$.
This result relies on the assumption, that
$\Rk(\bd{P}) = N-1$
and
$\bd{n}\{1\} \ne 0$.
We find that $\mathfrak{L}$ is a solvable Lie algebra
with an Abelian derived subalgebra;
if $\bd{P}$ is strictly upper-triangular,
then the Lie algebra~$\mathfrak{L}$ is nilpotent as well.

Relation (\ref{eq:adjmatfromt3}) implies
that all matrices~$\bd{A}_k$
share the same null space.
Furthermore, it is worthwhile to note that
\begin{eqnarray}
\label{eq:prodtjtimestk}
\bd{T}_j \cdot \bd{T}_k = \bd{n}\{j\}\;\bd{T}_k
\quad;
\end{eqnarray}
left-multiplication by $\bd{P}$ yields
\begin{eqnarray}
\bd{A}_j \cdot \bd{T}_k = \bd{n}\{j\}\;\bd{A}_k
\quad.
\end{eqnarray}
I.e.\ the matrix~$\bd{T}_k$
not only generates
the adjoint representation matrix~$\bd{A}_k$
from the a-priori matrix~$\bd{P}$,
it also may by used to transform
any adjoint representation matrix
into~$\bd{A}_k$.

\printbibliography

\appendix
\section{Computer implementation}
\label{sect:appendix}

An implementation of
the algorithm described in section~\ref{sect:algorithm}
using the Octave programming language
\parencite{eaton:2008}
is reproduced below.

\noindent
\rule{\textwidth}{0.4pt}

\renewcommand{\baselinestretch}{1.0}
\begin{verbatim}
function [adjRep, structConst] = randsolvableliealg( nofDim, allowCplx)



  if nargin < 1 || isempty( nofDim), nofDim = 3; end
  if nargin < 2 || isempty( allowCplx),
    allowCplx = false;
  end
  aPriori = randn( nofDim, nofDim-1);
  if allowCplx,
    aPriori = aPriori + 1i * randn( nofDim, nofDim-1);
  end
  % note:
  %  aPriori = aPriori - tril( aPriori);
  % generates a random nilpotent Lie algebra
  if rank( aPriori) ~= nofDim-1,
    error( [upper( mfilename) ...
      ': null space of parameter matrix not one-dimensional.'])
  end
  pMat          = zeros( nofDim, nofDim);
  % first column of 'pMat' is zero and ...
  pMat(:,2:end) = aPriori;
  % ... thus null space of matrix 'aPriori' is one-dimensional
  nVct    = transpose( null( transpose( aPriori)));
  assert( size( nVct, 1) == 1)
  unitMat = eye( nofDim);
  adjRep  = zeros( nofDim, nofDim, nofDim) * NaN;
  for k = 1:nofDim,
    adjRep(:,:,k) = pMat * (nVct(k) * unitMat - unitMat(:,k) * nVct);
  end
  if nargout > 1,
    structConst = permute( adjRep, [3, 2, 1]);
  end
end
\end{verbatim}

\vspace*{-4mm}
\noindent
\rule{\textwidth}{0.4pt}
\vspace*{-1mm}
\renewcommand{\baselinestretch}{1.2}

\noindent
On a standard personal computer
(2~GHz Intel(R) Core(TM)2 Duo processor, 4~GByte memory)
running a Linux operating system
with GNU Octave (version~3.8.1)
it takes about 0.3~seconds
to generate a solvable Lie algebra with dimension~$N=100$.
The solution for $N=500$ requires about 40~seconds.

\end{document}